\documentclass[11pt]{article}
\usepackage[usenames]{color}
\usepackage[colorlinks=true,
linkcolor=webgreen,
filecolor=webbrown,
citecolor=webgreen]{hyperref}

\definecolor{webgreen}{rgb}{0,.5,0}
\definecolor{webbrown}{rgb}{.6,0,0}

\usepackage{amssymb,psfig,epsfig,latexsym,graphicx,here}

\newtheorem{lemma}{Lemma}
\newtheorem{theorem}{Theorem}
\newtheorem{coro}{Corollary}
\newtheorem{conj}{Conjecture}

\def\binom#1#2{{#1}\choose{#2}}

\newcommand{\eqn}[1]{(\ref{#1})}
\newcommand{\Stop}{\mbox{Stop}}
\newcommand{\Prob}{\mbox{Prob}}
\newcommand{\al}{\alpha}
\newcommand{\be}{\beta}

\newcommand{\de}{\delta}
\newcommand{\eps}{\varepsilon}
\newcommand{\Om}{\Omega}
\newcommand{\sset}{\subseteq}
\newcommand{\bsq}{{\vrule height .9ex width .8ex depth -.1ex }}

\newcommand{\ZZ}{{\mathbb Z}}
\newcommand{\QQ}{{\mathbb Q}}

\newcommand{\sP}{{\mathcal P}}
\newcommand{\eeq}{\end{equation}}
\newcommand{\beql}[1]{\begin{equation}\label{#1}}

\newcommand{\leftBra}{\{ \hspace*{-.10in} \{ }
\newcommand{\rightBra}{\} \hspace*{-.10in} \} }

\makeatletter
\def\@sect#1#2#3#4#5#6[#7]#8{\ifnum #2>\c@secnumdepth
     \def\@svsec{}\else
     \refstepcounter{#1}\edef\@svsec{\csname the#1\endcsname.\hskip .75em }\fi
     \@tempskipa #5\relax
      \ifdim \@tempskipa>\z@
        \begingroup #6\relax
          \@hangfrom{\hskip #3\relax\@svsec}{\interlinepenalty \@M #8\par}%
        \endgroup
       \csname #1mark\endcsname{#7}\addcontentsline
         {toc}{#1}{\ifnum #2>\c@secnumdepth \else
                      \protect\numberline{\csname the#1\endcsname}\fi
                    #7}\else
        \def\@svsechd{#6\hskip #3\@svsec #8\csname #1mark\endcsname
                      {#7}\addcontentsline
                           {toc}{#1}{\ifnum #2>\c@secnumdepth \else
                             \protect\numberline{\csname the#1\endcsname}\fi
                       #7}}\fi
     \@xsect{#5}}
\def\@begintheorem#1#2{\it \trivlist \item[\hskip \labelsep{\bf #1\ #2.}]}

\def\section{\@startsection {section}{1}{\z@}{-3.5ex plus -1ex minus
 -.2ex}{2.3ex plus .2ex}{\normalsize\bf}}
\makeatother
\makeatletter
\def\subsection{\@startsection {subsection}{1}{\z@}{-3.5ex plus -1ex minus
 -.2ex}{2.3ex plus .2ex}{\normalsize\bf}}

\makeatother

\begin{document}
\begin{center}
{\large\bf Approximate Squaring} \\
\vspace*{+.2in}
{\em J. C. Lagarias} and {\em N. J. A. Sloane} \smallskip \\
Information Sciences Research Center \\
AT\&T Shannon Lab \\
Florham Park, NJ 07932--0971 \medskip \\
Email addresses: \href{mailto:jcl@research.att.com}{{\tt jcl@research.att.com}}, \href{mailto:njas@research.att.com}{{\tt njas@research.att.com}}
\bigskip \\
Aug. 16, 2003; revised Nov. 9, 2003 \bigskip
\end{center}

\begin{center}
{\bf Abstract} \\
\end{center}
We study the ``approximate squaring'' map
$f(x) := x \lceil x \rceil$ and its behavior when iterated.
We conjecture that if $f$ is repeatedly applied to a 
rational number $r = l/d > 1$ then
eventually an integer will be reached.
We prove this when $d=2$, and
provide evidence that it is true in general by giving an upper 
bound on the density of the ``exceptional set''
of numbers which fail to reach an integer.
We give similar results for a $p$-adic analogue
of $f$, when the exceptional set is nonempty,
and for iterating the ``approximate multiplication'' map
$f_r(x) := r \lceil x \rceil$, where $r$ is a fixed
rational number. We briefly discuss what happens when ``ceiling''
is replaced by ``floor'' in the definitions. \\

\vspace{0.7\baselineskip}
AMS 2000 Classification: Primary 26A18; Secondary: 11B83, 11K31, 11Y99 \\

%
%
%
%

\section{Introduction}\label{Sec1}
In this paper we study the ``approximate squaring'' map
$f : \QQ \rightarrow \QQ$ given by
\beql{Eq1.1}
f(x) := x \lceil x \rceil ~,
\eeq
and consider its behavior when iterated.
Although there  is an extensive literature
on iterated maps (see for example [Collet and Eckmann 1980],
[Beardon 1991], [Lagarias 1992]), 
including 
the study of various  first-order recurrences involving the 
ceiling function  ([Eisele and Hadeler 1990], [Graham and Yan 1999]),
the approximate squaring map  seems not to have 
been treated before, and has some interesting 
features. 

The function $f$ behaves qualitatively
like iterating the rational function 
$R(x) = x^2$.  Indeed, all points $|x| \le 1$ have a 
bounded orbit under $f(x)$, while all points
$|x| > 1$ have unbounded orbits and diverge
to $\infty$, just as they do when $R(x)$ is iterated.
However, $f(x)$ has the additional feature that it is
discontinuous at integer points.
It follows that the $n$-th iterate $f^{(n)}$ is discontinuous at
a certain set of rational points, namely, 
those points $x$ where $f^{(n)}(x)$ is an integer.

It is therefore natural to ask: if we start
with a rational number $r$ and iterate $f$,
will we always eventually reach an integer?
This question is the subject of our paper.

Numerical experiments suggest that the answer 
to our question is ``Yes'', although 
it may take many steps, and consequently involve
some very large numbers.

For example, starting at $r = \frac{3}{2}$, 
$f(r) = \frac{3}{2} \cdot 2 = 3$, reaching an integer
in one step; and 
starting at $r = \frac{8}{7}$ we get
$f(r) = \frac{16}{7}$,
$f^{(2)}(r) = \frac{48}{7}$,
$f^{(3)}(r) = 48$, taking three steps.
On the other hand, starting at $r = \frac{6}{5}$, we find
$$
f(r) = \frac{12}{5},~
f^{(2)}(r) = \frac{36}{5},~
f^{(3)}(r) = \frac{288}{5},~
f^{(4)}(r) = \frac{16704}{5},~
f^{(5)}(r) = \frac{55808064}{5},~
$$
$$
f^{(6)}(r) = \frac{622908012647232}{5}, ~
f^{(7)}(r) = \frac{77602878444025201997703040704}{5}, ~
\ldots,
$$
and we do not reach an integer until $f^{(18)}(r)$,
which is a number with $57735$ digits.

We note that for any rational starting point $r$, since $\lceil x \rceil$
is an integer, the denominators $d_j$ of the iterates $f^{(j)}(r)$ 
must form a nonincreasing sequence with $d_{j+1}$ dividing
$d_j$. For $0 < r \le 1$, $\lceil r \rceil = 1$ and $f(r) = r$,
so there the denominator is fixed.  
For $-1 < r \le 0$, $f(r) = 0$,
and for $r \le -1$, $f(r) \ge 1$.
So it is sufficient to restrict our attention to the
case of rationals $r > 1$.

We make the following conjecture:

\begin{conj}\label{Con1}
For each rational $r \in \QQ$ with $r > 1$,
there is an integer $m \ge 0$ such that
$f^{(j)}(r)$ is an integer for all $j \ge m$.
\end{conj}

We establish the conjecture in the
special case when the denominator is 2, where a complete analysis
is possible. This is done in Section \ref{Sec2}.

In Section \ref{Sec3} we consider the
case of rational starting values $r$
with a fixed  denominator $d \ge 3$. We show that
the set of  starting
values that reach an integer after exactly $j$ steps has
a limiting density, and that the set of starting
values that never reach an integer has density zero.
More precisely, at most a sparse subset of the rationals
$\{ r = \frac{l}{d} : d < l \le x \}$
fail to become integers under iteration, in
the sense that the cardinality of this subset is bounded above
by $C(d, \epsilon)x^{1 - \alpha_d + \epsilon}$ for a certain positive 
constant $\alpha_d$ and any positive $\epsilon$, where
$C(d, \epsilon)$ is a positive constant depending only on $d$ and $\epsilon$.
Showing that this ``exceptional set'' of
starting values which fail to reach integers is in fact empty (or even finite)
appears to be a difficult problem, for reasons indicated below.
We also show that the set of  starting
values that reach an integer after exactly $j$ steps has
a limiting density.

In Section \ref{Sec4} we consider a $p$-adic analogue of the
approximate squaring map. In this case we show that there
is a nonempty exceptional set of elements in $\frac{1}{p^k} \ZZ_p$
which under iteration never ``escape'' to the
smaller invariant set $ \frac{1}{p^{k-1}} \ZZ_p$.
This set has Hausdorff dimension
exactly $1 - \alpha_{p^k}$, where
$\alpha_{p^k}$ is the same constant that appeared in Section \ref{Sec3}.
The existence of this 
exceptional set is one reason why it may be a difficult
problem to obtain better upper bounds 
on the cardinality
of the exceptional set in Section \ref{Sec3}.

In Section \ref{Sec5} we study 
similar questions concerning 
the ``approximate multiplication'' map 
\beql{Eq1.2}
f_r(x) := r \lceil x \rceil ~,
\eeq
where $r$ is a fixed rational number.
For $r = \frac{1}{b}$, this map is a special case of the
map $x \mapsto a + \lceil \frac{x}{b} \rceil$,
where $a, b \in \ZZ$, $b \ge 2$, studied 
by P. Eisele and R. P. Hadeler [Eisele and Hadeler 1990].
Recently, J. S. Tanton [Tanton 2002],
together with Charles Adler, 
formulated  a game-theoretic problem ``Survivor'', 
and noted that its analysis leads to the study of 
the sequence of rational numbers
$$
a_{0} = r, ~ a_n = r \lceil a_{n-1} \rceil \mbox{~for~} n \ge 1 ~,
$$
for $r > 1$, which is the trajectory of $r$ under the map $f_r(x)$.
He raised the question, "Must some $a_n$ be an integer?",
and conjectured that the answer is ``Yes''. 
This question  differs from the case
of the approximate squaring map in that the denominators
of successive iterates, though bounded by the denominator of $r$,
may increase or decrease. We note that the long-term
dynamics of iterating this map differs according to
whether $|r| < 1$, $|r|=1$ or $|r|> 1$, with the
case $r > 1$ being most analogous to the approximate squaring
map.

The approximate multiplication maps have some 
resemblance to the map occurring in the $3x+1$ problem.
Setting  $r= \frac{l}{d}$, 
we observe that $f_r(x)$ maps the domain
$\frac{1}{d}\ZZ$ into itself, and on this domain is 
conjugate to the map $g_r: \ZZ \to \ZZ$ given by
\[
g_{r}(n) = \left\{
 \begin{array}{lll}
   \frac{1}{d}n &  \mbox{~if~} n  \equiv 0~(\bmod~ d) \, , \\
   \frac{1}{d}(n + l(d-b))&  \mbox{~if~} n   \equiv b~(\bmod~ d), 1 \le b \le d-1\, ,
 \end{array}
\right.
\]
(see \eqn{Eq5.1aa} of Section \ref{Sec5}).
In terms of the conjugated map the  question we consider becomes  
whether for most starting values
some iterate of  $g_r$ is an integer divisible by $d$.
For $r= \frac{3}{2}$ the conjugated map is
\[
g_{3/2}(n) = \left\{
 \begin{array}{ll}
   \frac{3}{2}n &  \mbox{~if~ $n \equiv 0~(\bmod~ 2)$} \, , \\
   \frac{3}{2}n + \frac{3}{2}&  \mbox{~if~ $n \equiv 1~(\bmod~ 2)$}\,.
 \end{array}
\right.
\]
This is similar in form to the $3x+1$ function
\[
T(n) = \left\{
 \begin{array}{ll}
   \frac{1}{2}n &  \mbox{~if~ $n \equiv 0~(\bmod~ 2)$} \, , \\
   \frac{3}{2}n + \frac{1}{2}&  \mbox{~if~ $n \equiv 1~(\bmod~ 2)$}\,.
 \end{array}
\right.
\]
as given  in  [Lagarias 1985] and [Wirsching 1998],
although the long-term dynamics of $g_{3/2}$
and $T$ are different.

We formulate a conjecture for approximate multiplication
maps analogous to the one above for the approximate squaring map.
We define the 
exceptional set $E(r)$ for the 
map $f_r$ to be
\beql{Eq1.aa}
E(r) := \{ n:~ n \in \ZZ~~ \mbox{and~ no~ iterate}~~ 
f_r^{(j)}(n) \in \ZZ \mbox{~for~} j \ge 1 \} ~.
\eeq
Then we have:

\begin{conj}\label{Con2}
For each nonintegral rational $r \in \QQ$ with $|r| > 1$,
the exceptional set $E(r)$ for the approximate
multiplication map $f_r$ is finite.
\end{conj}
The ``expanding map'' hypothesis $|r| > 1$ is necessary
in the statement of
this conjecture, for  the conclusion 
fails for $r= \frac{1}{d}$ with $d \ge 3$,
as remarked below. In parallel to the
results for the approximate squaring map we prove
Conjecture~\ref{Con2} 
for rational $r$ having denominator $2$; it 
remains open for all rationals 
with $|r| > 1$ having  denominator $d \ge 3$.
This conjecture may also be difficult, as
indicated 
in Section \ref{Sec5}
by an analogy with the problem of showing 
there exist no Mahler $Z$-numbers,
a notorious problem connected with 
powers of $\frac{3}{2}$ ([Mahler 1968], [Flatto, 1991];
also [Choquet 1980],  [Lagarias 1985]).

Our analysis in  Section \ref{Sec5} applies more generally to
the family $\sP_r$ of maps $h_r: \ZZ \to \ZZ$
having the form 
\beql{Eq1.4}
h_r(n) = \frac{1}{d}( l n + l_b) \mbox{~when~} n \equiv b~ (\bmod~d) ~,
\eeq
where the integers $l_b$ satisfy $l_b \equiv - l b~~(\bmod~d)$.
We obtain for all functions $h_r \in \sP_r$
an explicit  upper bound on the cardinality of the 
exceptional set $E(h_r, x)$ consisting of all integers
$|n| \le x$ that do not have some iterate that is a multiple
of $d$. 
We show that for all rationals $r$, 
$$
\#E(h_r, x) \le 4d x^{\beta_d},  
$$
with $\beta_d = \frac{\log(d-1)}{\log d}.$ 
We complement this upper bound with a result showing  that
it is of the
correct order of magnitude (to within a multiplicative
constant) for certain
values of $r$ lying in $0 < r < 1$.
This is the case
for the function $g_r$  for  $r = \frac{1}{d}$ with $d \ge 3$.
It follows that Conjecture~\ref{Con2}
cannot hold for these values of $r$. 

Section \ref{Sec6} gives  some numerical results related to these questions.

As a final point we  
observe that one can also consider similar functions 
obtained by replacing
the ``ceiling'' function with  the ``floor'' function.
Exactly the same set of questions can be asked  for such
functions. Our methods  carry over to such functions, and 
there appear to be analogues of all our results and conjectures.
For example:

\begin{conj}
Let $F(x) := x \lfloor x \rfloor$.
For each rational $r \in \QQ$ with $r \ge 2$,
there is an integer $m \ge 0$ such that
$F^{(j)}(r)$ is an integer for all $j \ge m$.
\end{conj}

We will not say any more about the ``floor'' function
versions, apart 
from one result at the end of Section \ref{Sec5}.

\paragraph{Notation:} 
$\lceil ~ \rceil$ denotes the ceiling function,
$\lfloor ~ \rfloor$ the floor function, and
$\leftBra ~ \rightBra$ the fractional part.
For a prime $p$, $|~|_p$ denotes the $p$-adic
valuation. (If $r \in \QQ \, , \, r = p^a \, \frac{b}{c}
\mbox{~with~} a, b, c \in \ZZ, c \neq 0, \gcd(p,b) = \gcd(p,c) = 1,
\mbox{~then~} |r|_p = a$.)
$\QQ _p$ and $\ZZ _p$ denote the $p$-adic rationals and integers,
respectively.
For integers $r, s, i$, $r|s$ means $r$ divides $s$,
and $r^i ||s$ means $r^i$ divides $s$ but $r^{i+1}$ does not.

%
%
%
%
\section{Denominator $2$}\label{Sec2}
In this section we investigate the case when the starting value $r$ has
denominator $2$. Here we are 
able to give a complete analysis.
The following table shows what happens for the first few
values of $r$.
It gives the initial term, the number of steps to reach an integer,
and the integer that is reached.

$$
\begin{array}{rcccccccccc}
\mbox{start}: & \frac{3}{2} & \frac{5}{2} & \frac{7}{2} & \frac{9}{2} & \frac{11}{2} & \frac{13}{2} & \frac{15}{2} & \frac{17}{2} & \frac{19}{2} & \cdots \\
\mbox{steps}: & 1 & 2 & 1 & 3 & 1 & 2 & 1 & 4 & 1 & \cdots \\
\mbox{reaches}: & 3 & 60 & 14 & 268065 & 33 & 2093 & 60 & 1204154941925628 & 95  & \cdots
\end{array}
$$

The number of steps appears to match sequence \htmladdnormallink{A001511}{http://www.research.att.com/cgi-bin/access.cgi/as/njas/sequences/eisA.cgi?Anum=A001511}
in [Sloane 1995--2003] (and the numbers reached now form sequence
\htmladdnormallink{A081853}{http://www.research.att.com/cgi-bin/access.cgi/as/njas/sequences/eisA.cgi?Anum=A081853} in that database). Indeed, we have:

\begin{theorem}\label{Th0}
Let $r = \frac{2l+1}{2}$, with $l \ge 1$. Then $f^{(m)}(r)$
reaches an integer for the first time when $m = |l|_2 + 1$.
\end{theorem}

\paragraph{Proof.}
Note that if $x \in \QQ$ has denominator $2$ and
is not an integer then $\lceil x \rceil = x + \frac{1}{2}$. 

We use induction on $|l|_2 = v$. If $v = 0$ then 
$l$ is odd, $\lceil r \rceil = r + \frac{1}{2} = l + 1$ is even,
and $r \lceil r \rceil$ has become an integer in one step,
as claimed.

Suppose $v \ge 1$, and
$$
l = 2^v + l_{v+1} 2^{v+1} + l_{v+2} 2^{v+2} + \cdots ~
$$
is the binary expansion of $l$,
where each $l_i = 0$ or $1$. Then
\begin{eqnarray}
r \lceil r \rceil & ~=~ & \left ( l + \frac{1}{2} \right ) ( l + 1) \nonumber \\
& ~=~ & \frac{1}{2} + \frac{l}{2} + l + l^2  \nonumber \\
& ~=~ & \frac{1}{2} + 2^{v-1} + (l_{v+1}+1) 2^v + 
(l_{v+1} + l_{v+2}) 2^{v+1} + \cdots + 2^{2v} + 
\cdots  \nonumber \\
& ~=~ & \frac{2 l' + 1}{2}  \nonumber
\end{eqnarray}
where
$$
l' = 2^{v-1} + (l_{v+1}+1) 2^{v} + (l_{v+1} + l_{v+2}) 2^{v+1} + \cdots + 2^{2v} + \cdots ~.
$$
and $|l'|_2 = v-1$. By the induction hypothesis, this will reach an integer in $v-1$ steps,
so we are done. ~~~$\bsq$

\paragraph{Remark.}
The numbers $2l+1$ for which $|l|_2 = v$ are precisely the numbers 
that are congruent to $2^{v+1} + 1 ~(\bmod~ 2^{v+2})$.
For example, if $v=0$, $2l+1 \in \{3, 7, 11, 15, \ldots\}$,
of the form $3$ (mod $4$); if $v=1$,
$2l+1 \in \{5, 13, 21, 29, \ldots\}$,
of the form 5 (mod 8); and so on.

\begin{coro}
Let $r = \frac{2l+1}{2}, l \ge 1, |l|_2 = v$.
Then the first integer value taken by $f^{(m)}(r)$ is
$$
\frac{1}{2} \, \theta ^{(v+1)}(2l+1) ~,
$$
where $ \theta(y) = y(y+1)/2$.
\end{coro}

\paragraph{Proof.}
This is now a straightforward calculation, again
using the fact that
if $x \in \QQ \setminus \ZZ$ has denominator 2
then $\lceil x \rceil = x + \frac{1}{2}$.  ~~~$\bsq$

For example, if $v=0$, and $r= (4k+3)/2 = y/2$ (say),
then in one step we reach the integer $\frac{1}{2} \theta(y) = y(y+1)/4$.
If $v=1$, and $r= (8k+5)/2 = y/2$,
then in two steps we reach the integer 
$$
\frac{1}{2} \theta(\theta(y)) = \frac{ y(y+1)(y^2+y+2)}{16} ~;
$$
if $v=2$, and $r= (16k+9)/2 = y/2$,
then in three steps we reach the integer 
$$
\frac{1}{2} \theta(\theta(\theta(y))) = 
\frac {y\left (y+1\right )\left ({y}^{2}+y+2\right )
\left ({y }^{2}-y+2\right )
\left ({y}^{2}+3\,y+4\right )}{256} ~,
$$
and so on.

%
%
%
%

\section{Denominator $d$}~\label{Sec3}

We now analyze the case of 
rationals with a general denominator $d$,
obtaining less complete results. 
The next theorem shows that most rationals will
eventually reach an integer.  In particular,
it gives an upper bound on the number  of 
such rationals below $x$ 
that never reach an integer. 
Given an integer  $d \ge 2$, and a bound $x \ge 1$ 
we study the ``exceptional set'' 
\beql{Eq3.1}
\mathcal{M}_d(x) ~:=~  \{ l ~:~ 1 \le l \le x,~ f^{(m)}( \frac{l}{d}) \notin \ZZ
\mbox{~for~each~} m \ge 1 \} ~,
\eeq
and let $M_d(x) = | \mathcal{M}_d(x) |$.
The finite set $[1, d-1] :=\{1, 2, \ldots, d-1\}$
is contained in $\mathcal{M}_d(x)$, and 
Conjecture \ref{Con1} asserts that $\mathcal{M}_d(x) =[1, d-1].$

\begin{theorem}\label{Th1}
Let $d \ge 2$ be a fixed integer.

(1) For each finite $j \ge 0$ 
the  set of positive integers $l$ 
for which  the approximate squaring map with
initial value  $r = \frac{l}{d}$ 
first iterates to an integer after exactly $j$ steps 
is a union of 
arithmetic progressions $(\bmod~ d^{j+1})$ restricted
to positive integers.

(2) There is a positive  exponent $\alpha_d$
such that  for each $\eps > 0$ and all $x > 1$,
\beql{Eq3.2}
M_d(x) \le C(d, \eps) \, x^{1 - \al_d + \eps} ~
\eeq
for a positive constant $C(d, \eps)$, with $\alpha_d$ given by
\beql{Eq3.3}
\al_d ~=~ \min_{ d' | d, \, d' > 1} \log_{d'} \left ( \frac {d'} { \phi (d')} \right ) ~,
\eeq
where $\phi$ is the Euler totient function.
In fact
\beql{Eq3.4}
\al_d = \min_{ p^j || d}
\frac{ \log \left (1 + \frac {1}{p-1} \right )}{ j \log p} ~.
\eeq
\end{theorem}

\paragraph{Note:} It follows immediately
from \eqn{Eq3.3} that $0 < \al_d \le 1$ and $\al_d = 1$ only
for $d = 2$. \\

To prove the 
theorem we  shall first  prove a preliminary lemma. We 
need some notation concerning the pattern of denominators
in a sequence of iterates.
We write $\frac{l}{d} = \frac{l_0}{d_0}$
with $\gcd (l_0, d_0) = 1$, and set
$$
f^{(j)} \left ( \frac{l}{d} \right ) = \frac{l_j}{d_j} ~, ~ j \ge 1 ~,
$$
where $\gcd(l_j, d_j) = 1$ and $d_j | d$ with $d_j \ge 1$.
The pair $(l,d)$ determines the sequence $d_0, d_1, d_2, \ldots$.

We consider the first $m$ iterates for a given input $r= \frac{l}{d}$,
and let
\beql{Eq3.9a}
{\mathcal Y(l, d, m)} := (d_0, d_1, ..., d_m)
\eeq
denote the vector of  denominators of the first $m$ iterates.
We call this a {\em chain} of length $m+1$.
This information in ${\mathcal Y}$ can be described in another way. 
We set  $d_{-1}=d$ 
and define $r_{j}$ by
$$
r_j = \frac {d_{j-1}}{d_j}
$$
for $j = 0, 1, 2, \ldots, m\,$. We call $j$ a
{\em break-point} if $r_j > 1$. 
Let the break-points up to the $m$-th iterate
be $0 \le j_1 <  j_2 <  \ldots <  j_t \le m$.
The data ${\mathcal Y(l, d, m)}$
is completely determined by knowledge of the pairs
$(j_1, r_{j_1}), ..., (j_t, r_{j_t})$ and $m$. 
Note that since  each break-point removes a
divisor of $d$ from the denominator, in any
chain of length $m$ there are at most $s$ break-points,
 where $s$ is the number of prime factors of $d$,
counted with multiplicity.

Now consider an arbitrary chain ${\mathcal Y}$, which
consists of a sequence $(d_0, d_1,..., d_m)$ with
$d_{j+1}$ dividing $d_j$ and $d_0$ dividing $d$.
A chain is {\em complete} if $d_m=1$, and is
{\em incomplete} otherwise. 
We will only  consider complete
chains with $d_{m-1} \ge 2$, i.e. with
$j_t=m$. 

The following lemma characterizes the set of positive integers $l$ 
that have a given associated chain ${\mathcal Y}$ of length $m+1$.

\begin{lemma}~\label{le31a}
Given a fixed chain ${\mathcal Y}= (d_0, d_1, ..., d_m)$ of length $m+1$,
the set of positive integers $l$ such that $r= \frac{l}{d}$ has
${\mathcal Y}(l, d, m) = {\mathcal Y}$ 
consists of a 
collection  of arithmetic progressions \\
$(\bmod~d_{-1} d_0 d_1 \cdots d_{m-1})$,
restricted to  the positive integers. There here are exactly
\beql{Eq3.12}
\phi(d_0)\phi(d_1) \phi(d_2) \cdots \phi(d_m)
\eeq 
such arithmetic progressions.
\end{lemma}

\paragraph{Proof.} 
We study the set of  $r=\frac{l}{d}$ having a fixed chain ${\mathcal Y}$
for its  initial $m$ iterates.

The set of allowable $l$ with $(l, d) = d_0$ consists of all integers
$r_1 l_0$ with $gcd(l_0, d_0)$=1. This consists of $\phi(d_0)$
arithmetic progressions  $(\bmod~ d)\equiv(\bmod~ d_{-1})$. 
 
We now write $\frac{l_k}{d_k}$ (for $k = 0, 1, \ldots$)
in a mixed-radix
expansion where the radices depend on $k$:
\beql{Eq3.13}
\frac{l_k}{d_k} = 
\frac{a_{-1}(k)}{d_k} + a_0(k) + 
\sum_{j=1}^{\infty} a_j(k) \prod_{l=0}^{j - 1} d_{k+l}~,
\eeq
in which the ``digits'' $a_j(k)$ satisfy
$$
0 < a_{-1}(k) < d_k \, ,
\mbox{~and~}
0 \le a_{j}(k) < d_{j+k} \mbox{~for~each~} j \ge 0 ~. 
$$
Here we set $d_{m+j} := d_m$ for all $j \ge 1$.
The sum on the right-hand side of \eqn{Eq3.13}
is actually a finite sum. By definition,
\begin{eqnarray}\label{Eq3.16}
\frac{l_{k+1}}{d_{k+1}} ~=~ f \left ( \frac{l_k}{d_k} \right ) & ~=~  &
\left ( \frac{a_{-1}(k)}{d_k} + a_0(k) +
\sum_{j=1}^{\infty} a_j(k) d_k d_{k+1} \cdots d_{k+j-1} \right ) \nonumber \\
& ~ & ~~ \times ~~ 
\left ( 1 +  a_0(k) +
\sum_{j=1}^{\infty} a_j(k) d_k d_{k+1} \cdots d_{k+j-1} \right ) ~. 
\end{eqnarray}

We use induction on $k \ge 0$ to simultaneously establish four properties
of this mixed-radix expansion:

(i) We have
\beql{Eq3.17}
r_{k+1} ~=~ \gcd ( a_0(k)+1, d_k) ~.
\eeq

(ii) We have  
\beql{Eq3.18}
a_{-1}(k+1) ~ \equiv ~ 
a_{-1}(k) \frac{a_0(k)+1}{r_{k+1}} ~ ( \bmod~ d_{k+1} ) ~. 
\eeq

(iii) We have 
\beql{Eq3.19}
\gcd ( a_{-1}(k+1), d_{k+1}) ~=~ 1 ~.
\eeq

(iv) For  $0 \le j \le m$,
\beql{Eq3.20}
a_j(k+1) ~ \equiv ~ a_{-1}(k) a_{j+1}(k) + 
G( a_{-1}(k), a_0(k), a_1(k), \ldots, a_{j}(k)) ~ (\bmod~ d_{j+k+1}) ~,
\eeq
for some  function $G( a_{-1}(k), a_0(k), a_1(k), \ldots, a_{j}(k))$.
This function includes all the necessary information about ``carries''
in the multiple-radix expansion. 

The base case $k=0$ is checked directly.
Since $d_0$ divides all terms in the sum
on the right-hand side of \eqn{Eq3.13},
the right-hand side of \eqn{Eq3.16} (when $k=0$) has a single term
$\frac{ a_{-1}(0) (a_{0}(0)+1)}{d_0}$ having a denominator,
and this term equals $\frac{a_{-1}(1)}{d_1} (\bmod~1)$.
Since $\gcd ( a_{-1}(0), d_0) = 1$
by hypothesis, we must have
$$
\gcd ( a_0(0)+1, d_0) = \frac{d_0}{d_1} = r_1~,
$$
which is (i) for $k=0$.
The term with a denominator in \eqn{Eq3.16} is then
$$
\frac {a_{-1}(0) \frac{a_0(0)+1}{r_1} }{d_1} ~,
$$
hence
\beql{Eq3.20a}
a_{-1}(1) \equiv a_{-1}(0) \left ( \frac{a_0(0)+1}{r_1} \right ) 
~ (\bmod~ d_1) ~,
\eeq
which is (ii).
Now $\gcd ( a_{-1}(1), d_1)$
divides $\gcd ( a_{-1}(0), d_1) \gcd( \frac{a_0(0)+1}{r_1}, d_1)$,
both terms of which are $1$, so (iii) follows.
Finally, to establish (iv)
when $k=0$,
we drop the terms involving $d_1 d_2 \cdots  d_{j+1}$
from \eqn{Eq3.16} and
observe that there is a  term
\beql{Eq3.21}
a_{-1}(0) a_{j+1}(0) d_1 d_2 \cdots d_j ~,
\eeq
while  all the other terms containing any
$a_{j+l}(0)$ for $l \ge 1$ are divisible by $d_1 d_2 \cdots d_{j+1}$.
(Note that $d_1 d_2 \cdots d_{j+1}$ divides $d_0 d_1 d_2 \cdots d_{j}$.)
This establishes that a congruence of the form \eqn{Eq3.20} 
holds for the
digit $a_j(1)$ in the expansion \eqn{Eq3.16}
for $\frac{l_1}{d_1}$, completing the proof of the base case.

The induction step for general $k$ follows using exactly the same reasoning.

Next, we will deduce from  (iv) 
that for each $k \ge 0$, 
\beql{Eq3.22}
a_0(k+1) ~ \equiv ~   a_{k+1}(0)  \prod_{l=0}^k a_{-1}(l) +
\widetilde{G}( a_{-1}(0), a_0(0), a_1(0), \ldots, a_{k}(0)) ~ (\bmod~ d_{k+1}) ~.
\eeq
for some function $\widetilde{G}$ depending on the indicated
variables. The important point about \eqn{Eq3.22} is that the dependence
on $a_{k+1}(0)$ is linear, even though the dependence on the
other initial terms $a_{-1}(0), a_0(0), \ldots, a_k(0)$ is
nonlinear.

To prove \eqn{Eq3.22}, we again use induction on $k$.
For $k=0$ the assertion is that
\beql{Eq3.24a}
a_0(1) ~ \equiv ~   a_{1}(0) a_{-1}(0) +
\widetilde{G}( a_{-1}(0), a_0(0)) ~ (\bmod~ d_{1}) ~,
\eeq
which is \eqn{Eq3.20} with $k=j=0$.
For $k=1$ we wish to show 
\beql{Eq3.24b}
a_0(2) ~ \equiv ~   a_{2}(0) a_{-1}(0) a_{-1}(1) +
\widetilde{G}( a_{-1}(0), a_0(0), a_1(0)) ~ (\bmod~ d_{2}) ~.
\eeq
Setting $k=1, j=0$ and $k=0, j=1$ in \eqn{Eq3.20} we obtain
$$
a_0(2) ~ \equiv ~   a_{-1}(1) a_{1}(1) +
\widetilde{G}( a_{-1}(1), a_0(1)) ~ (\bmod~ d_{2})
$$
and
$$
a_1(1) ~ \equiv ~   a_{-1}(0) a_{2}(0) +
\widetilde{G}( a_{-1}(0), a_0(0), a_1(0)) ~ (\bmod~ d_{2}) ~,
$$
hence
$$
a_0(2) ~ \equiv ~   a_{2}(0) a_{-1}(0) a_{-1}(1) +
a_{-1}(1)\widetilde{G}( a_{-1}(0), a_0(0), a_1(0)) +
\widetilde{G}( a_{-1}(1), a_0(1)) ~ (\bmod~ d_{2}) ~.
$$
However, from \eqn{Eq3.20a} and the fact that $0 < a_{-1}(1) < d_1$,
$a_{-1}(1)$ is uniquely determined by $a_{-1}(0)$ and $a_0(0)$.
Also the induction hypothesis allows us to use \eqn{Eq3.24a} to eliminate
$a_0(1)$. Equation \eqn{Eq3.24b} follows.
The case of general $k$ follows in the same way;
we leave the details to the reader.

We have already seen that
$$
\gcd ( a_{-1}(0), d_1) = 1 ~.
$$
From \eqn{Eq3.19} and \eqn{Eq3.18} ,
$$
\gcd ( a_{-1}(2), d_2) = 1 ~,
$$
$$
a_{-1}(2) ~ \equiv ~ a_{-1}(1) \frac{a_0(1)+1}{r_{2}} ~ ( \bmod~ d_{2} ) ~.
$$
Therefore $\gcd ( a_{-1}(1), d_2) = 1$ and so
$$
\gcd ( a_{-1}(0) a_{-1}(1), d_2) = 1 ~.
$$
Continuing in this way we obtain
\beql{Eq3.23}
\gcd (\prod_{l=0}^k a_{-1}(l), d_{k+1}) ~=~ 1 ~,
\eeq
for $k=0, 1, \ldots, m$.

Then \eqn{Eq3.22} shows
that the congruence class of $a_0(k+1) ~(\bmod~ d_{k+1})$ is
uniquely determined by the congruence class of $a_{k+1}(0) ~(\bmod~ d_{k+1})$,
once $a_{-1}(0), a_{0}(0), \ldots, a_k(0)$ are specified.

In particular for each $k \ge 0$
there are exactly $\phi(d_{k+1})$ congruence classes
of $a_{k}(0)~( \bmod~d_{k} )$ that give 
$$
\gcd ( \frac{a_0(k)+1}{r_{k+1}}, d_{k+1} ) ~=~ 1 ~,
$$
or in other words which give
$$
\gcd ( a_0(k)+1, d_{k} ) ~=~ r_{k+1} ~.
$$
At each iteration we impose one such condition.
Combining this with the $\phi(d_0)$ congruence classes
 $(\bmod~d_{-1})$
allowed initially, 
 \eqn{Eq3.12} follows. ~~~$\bsq$

\paragraph{Proof of Theorem~\ref{Th1}.}
Part (1) follows from Lemma~\ref{le31a}.
Indeed the  set of integers that iterate to an integer in
exactly $j$ steps are precisely those integers that
belong to some complete chain ${\mathcal Y}$ having $j= j_t$.
There are 
finitely many complete chains of length $m+1$ with $m= j_t = j$, and by
 Lemma~\ref{le31a} the set of positive $l$ such
that $\frac{l}{d}$ belongs to a given such chain 
is a union of 
arithmetic progressions $(\bmod~ d_{-1}d_0 d_1 \ldots d_{m-1})$.
Each such progression subdivides into 
a finite union of arithmetic progressions
$(\bmod~d^{m+1})$, since $ d_{-1} d_0 d_1 \ldots d_{m-1}$ divides $d^{m+1}.$

To establish part (2), let the cutoff value  $x$ be given.
We call a value $l \le x$ ``bad'' at size $x$ if
there is an iterate $m \ge 0$ such that
\beql{Eq3.5}
\prod_{j=0}^{m-1} d_j \le x < \prod_{j=0}^{m} d_j ~,
\eeq
with $d_m > 1$; here $m$ depends on $l, d$ and $x$.
The chain ${\mathcal Y}$ associated to
the data $(d_0, d_1, ..., d_m)$ is
necessarily an incomplete chain, and we say it is
{\em incomplete at size $x$} if \eqn{Eq3.5} holds. 
Conversely all elements $l \le x$ 
belonging to such an  incomplete chain
are ``bad''. We
let $N_d(x)$ denote the number of ``bad'' elements
at size $x$, and we clearly have
$$
M_d(x) \le N_d(x) ~.
$$
We will establish part (2) by showing the stronger result that
\beql{Eq3.7}
N_d(x) \le C(d, \eps) \, x^{1 - \al_d + \eps} ~.
\eeq

To bound $N_d(x)$ we first bound the number of distinct incomplete
chains at size $x$, and then bound  the maximal number
of elements below $x$ falling in any such chain.
It is easy to bound the  number of 
incomplete  chains ${\mathcal Y}$
at size $x$ as follows. Since each $d_j \ge 2$
we have $m \le \log_2~ x$.
There are ${\binom{m}{t}}$ choices for the break-points,
and the bound $t \le s$ 
then yields that there are at most
$$
\sum_{t=0}^{s} {\binom{\log_2 x}{t}}
\le s ( \log_2 x)^{s}
$$
break-point patterns. Since  
a nontrivial divisor $r_i$ of $d_{i-1}$ is removed at each break-point, we
conclude that the number of distinct chains is at most
$ d! s (\log_2 x)^{s}. $
For $d$ (and hence $s$) fixed  this gives for any fixed
$\epsilon > 0$ and $x \ge 2$ that
\beql{Eq3.11b}
 \#({\mathcal Y})  \le C_1(d, \eps) x^\eps.
\eeq
for some positive constant $ C_1(d, \eps)$.

Now let 
 $N_d(\mathcal Y, x)$ count the number of bad $l$ at size $x$ having
a given chain ${\mathcal Y}= (d_0, d_1, ..., d_m)$ 
that is incomplete at size $x$.
In view of \eqn{Eq3.11b}
it  suffices to
prove an upper  bound of the same form \eqn{Eq3.7} for each such
chain.
So consider a fixed such  chain ${\mathcal Y}$.
Now  \eqn{Eq3.5} holds, and
shows that each arithmetic progression 
$(\bmod ~d_0 d_1 \cdots d_{m-1}$) contains at most $d_m$ elements
below $x$.
From \eqn{Eq3.12}, $N_d( {\mathcal Y}, x)$, the number
of bad elements $l$ at size $x$ with chain $\mathcal Y$ ,
satisfies
$$
N_d( {\mathcal Y}, x) \le d_m \phi(d_1) \phi(d_2) \cdots \phi(d_m) ~.
$$
Now $x \ge d_0 d_1 \cdots d_{m-1} \ge d_1 d_2 \cdots d_{m-1}$,
hence
$$
\frac{N_d( {\mathcal Y}, x)}{x} \le 
d_m \prod_{j=1}^{m} \frac{\phi(d_j)}{d_j} ~,
$$
and therefore
\begin{eqnarray}\label{Eq3.27}
\frac{N_d( {\mathcal Y}, x)}{x} & ~\le~ &
d_m \, e^{ \sum_{j=1}^{m} \log \frac{\phi(d_j)}{d_j} } \nonumber \\
& ~\le~ &
d \, x^{ \frac { \sum_{j=0}^{m} \log \frac{\phi(d_j)}{d_j} }  { \log x } } ~.
\end{eqnarray}
Now $x \le d_0 d_1 \cdots d_m$, so that
\begin{eqnarray}\label{Eq3.28}
\frac
{ \sum_{j=0}^{m} \left | \log \frac { \phi(d_j) }{d_j} \right | }
{ \log x }
& ~\ge~ &
\frac
{ \sum_{j=0}^{m} \left | \log \frac { \phi(d_j) }{d_j} \right | }
{ \sum_{j=0}^{m} \log d_j } \nonumber \\
& ~\ge~ &
\min_{0 \le j \le m} ~
\frac
{ \left | \log \frac { \phi(d_j) }{d_j} \right | }
{ \log d_j } \nonumber \\
& ~\ge~ &
\min_{0 \le j \le m} ~
\left (
\log_{d_j}  \left ( \frac {d_j}{\phi(d_j)} \right )
\right )
\nonumber \\
& ~\ge~ & \al_d ~, \nonumber
\end{eqnarray}
where we used the definition \eqn{Eq3.3}.
Substituting this in \eqn{Eq3.27} yields
$$
\frac{N_d( {\mathcal Y}, x)}{x} \le 
d \, x^{-\al_d} ~,
$$
since $\frac{ \phi(d_j)}{d_j} < 1$ gives
$ \log \frac { \phi(d_j)}{d_j} < 0$,
and so
$$
N_d( {\mathcal Y}, x) \le 
d \, x^{1-\al_d} ~,
$$
as required.
Combined with \eqn{Eq3.11b}, this gives \eqn{Eq3.7}, hence \eqn{Eq3.2}.

Finally, we establish the equivalence of \eqn{Eq3.3}
and \eqn{Eq3.4}. Since
$\phi$ is multiplicative, for a general $d' > 1$ we have
\begin{eqnarray}\label{Eq3.31}
\log_{d'} \left ( \frac { d'}{\phi(d')} \right )
& ~=~ &
 \frac{\sum_{p^j || d'} \log \left ( \frac { p^j}{\phi(p^j)} \right )}
{\sum_{p^j || d} \log p^j}  \nonumber \\
& ~\ge~ &
\min_{p^j || d'}
\left \{ 
\log_{p^j} \left ( 
\frac { p^j}{\phi(p^j)} 
\right ) \right \} ~.
\end{eqnarray}
Thus the minimum in \eqn{Eq3.3} is attained
when $d'$ is a prime power. Now
\beql{Eq329a}
\log_{p^j} \left ( \frac { p^j}{\phi(p^j)} \right )
 ~=~ 
\frac
{\log \left ( 
\frac { 1}{1 - \frac{1}{p} }
\right ) 
}
{ \log~ p^j}
 ~=~ 
\frac
{\log \left (  1 + \frac { 1}{p-1}\right ) }
{ j \log p} 
\eeq
is minimized by making $j$ as large as possible,
so we obtain the formula \eqn{Eq3.4}.  ~~~$\bsq$

The upper bound in \eqn{Eq3.7} in Theorem~\ref{Th1}
has an exponent that is essentially best possible. Indeed, 
if $l, d, x, m$ are such that
$$
d_1 = d_2 = \cdots = d_m = d' ~,
$$
where $d'$ is the value that minimizes \eqn{Eq3.3}, 
so that $t = 1$ and the chain is simply
${\mathcal Y}= \{(1, d')\}$,
we have
$$
N_d( {\mathcal Y}, x) \ge
{C}_d^{'} \, x^{1-\al_d} ~,
$$
for some positive constant ${C}_d^{'}$.

We next give a probabilistic
interpretation of  Theorem~\ref{Th1},
following a suggestion of Daniel Berend.
Given a denominator $d$ and
an integer cutoff value $x$,
let $X$ be a random integer chosen 
uniformly on the interval  $[1, x]$,
and let $s = X/d$. We iterate $f$ starting at $s$,
and let $\Stop(X)$ be the random
variable giving the number
of iterates needed to reach an integer.
$\Stop(X)$ takes nonnegative integer values (including
$0$ and $\infty$). There is an  associated probability 
distribution
$$
\{ \Prob _x [\Stop = j]: j= 0, 1, \ldots , \infty \}
$$
which depends on the interval $[1,x]$.

\begin{theorem}~\label{Th3.2}
For each fixed  denominator $d$, there exists a 
limiting distribution as $x \to \infty$ of the stopping times
of the approximate squaring map $f$ on rationals $\frac{l}{d}$.
More precisely, for  each finite $j$ the limit 
\beql{Eq330}
\Prob _{\infty}[ \Stop(X) = j] := \lim_{x \to \infty} \Prob _x [\Stop(X)=j] 
\eeq
exists and is a rational number with denominator dividing $d^{j+1}$,
and for  $j= \infty$ we have  
\beql{Eq331}
 \Prob _{\infty}[ \Stop(X) = \infty] := 
 \lim_{x \to \infty} \Prob [\Stop(x) =  \infty] = 0~. 
\eeq
In the case when $d=p$ is  a prime the limiting distribution is 
a geometric distribution:
\beql{Eq332}
\Prob _{\infty}[ \Stop(X) = j] = p(1-p)^j ~.
\eeq
\end{theorem}

\paragraph{Proof.}
The assertion \eqn{Eq330} follows 
directly from part (1) of Theorem~\ref{Th1}, which shows
that the set with stopping time $j$ is a union of
complete arithmetic progressions $(\bmod~d^{j+1})$.  
The assertion \eqn{Eq331} follows from part (2)
of Theorem~\ref{Th1}.
The computation of the density \eqn{Eq332} follows from the
fact that for $d=p$ and each $j \ge 0$
there is a unique complete chain
${\mathcal Y}$ with stopping time $j$, having the single 
break-point  $(j, 1)$, giving one arithmetic progression $(\bmod~p^{j+1}).$
  ~~~$\bsq$

\paragraph{Remarks.}

(1) This stopping time formulation is similar
to that arising in the $3x+1$ problem
and studied in [Terras 1976, 1979].

(2) In principle the limiting distribution
for any specific composite $d$
is computable from the proof of Theorem~\ref{Th1}.

(3) For all $d$ we have
$
\Prob _{\infty}[\Stop(X) =0] = 1/d.
$

%
%
%
%

\section{$p$-Adic iteration}~\label{Sec4}

We now consider  an  approximate squaring  map 
defined on the $p$-adic numbers analogous to the
approximate squaring map on $\QQ$,
and study the question
of whether some iterate will eventually
become a $p$-adic integer. 
We show that now there is always a nonempty exceptional set
of $p$-adic numbers which never become $p$-adic integers.

Let $p$ be a prime $\ge 2$ and
let $\QQ_p$ denote the $p$-adic numbers,
with typical element $\al = \sum_{j=-k}^{\infty} a_j p^j$, 
where $k \in \ZZ$ and the $a_j$ satisfy $0 \le a_j \le p-1$.
The $p$-adic integral part of $\al$
is given by the function
$F_p: \QQ_p \rightarrow \QQ_p$, where
$$
F_p(\al) := \sum_{j= \max \{ 0, k \} }^{\infty} a_j p^j ~,
$$
while the $p$-adic fractional part (or ``principal part'') of $\al$ is
$$
P_p(\al) := \sum_{j <0} a_j p^j ~,
$$
which is a finite sum (possibly empty); thus 
$\al = P_p(\al) + F_p(\al).$
We investigate the  function
$$
f_p(\al) := \al ( F_p(\al) + 1) ~.
$$
which is a $p$-adic analogue of the approximate squaring map
defined in \eqn{Eq1.1}.
If we regard the rationals $\QQ$ as embedded in $\QQ_p$, then
for nonintegral $r$ in the subring
$\ZZ[\frac{1}{p}] \sset \QQ \sset \QQ_p$, we have  
$f_p(r) = r \lceil r \rceil,$
so the iterates there agree
with the approximate squaring map.

For each $k \ge 0$ the set $\frac{1}{p^k} \, \ZZ_p$
is invariant under the action of $f_p(x)$, and
$$
\ZZ \sset \frac{1}{p} \, \ZZ \sset \frac{1}{p^2} \, \ZZ \sset \cdots ~.
$$
We define the exceptional set
\beql{Eq4.5}
\Om_k(p) := \left \{
\al \in \frac{1}{p^k} \, \ZZ_p ~:~ 
f^{(m)}(\al) \notin \frac{1}{p^{k-1}} \, \ZZ_p \, \mbox{~for ~each~} ~ m \ge 1 \right \} ~,
\eeq
for $k \ge 1$.
The set $\Om_k(p)$ is an analogue
of the
exceptional sets $\mathcal{M}_d(x)$ studied in the last
section, corresponding to the denominator $d= p^k$.
(Note that  $\Om_k(p) \cap \frac{1}{p^k}\ZZ_{>0}$ is contained 
in  $\frac{1}{p^k}\mathcal{M}_{p^k}(\infty)$, as defined in \eqn{Eq3.1}).
We will show that
these sets are nonempty, determine their Hausdorff
dimension $s$, and get upper and lower bounds on their
Hausdorff $s$-dimensional measure.

The $s$-dimensional $p$-adic Hausdorff measure ${\mathcal H}_{p}^{s} ( \Om)$
of a closed set $\Om$ in $\QQ_p$ is defined by the general prescription in
[Falconer 1990, Chapter 2] or [Federer 1969, Section $\S 2.10$].
Here $0 < s \le 1$. The {\em diameter} of a measurable
set $S \sset \QQ_p$ is 
$$
|S| := \sup \{ |\al - \be|_p ~:~ \al, \, \be \in S \} ~,
$$
and its {\em $p$-adic measure} $\mu_p(S)$ is Haar measure with the
normalization $\mu_p(\ZZ_p) = 1$.
A $p$-adic ball
$$
B(\al; p^l) := \{ \be \in \QQ_p ~:~ | \be - \al|_p \le p^l \}
$$
is both closed and open, and has the property that its diameter
equals its measure:
$$
| B(\al; p^l) | = \mu ( B(\al; p^l) ) = \frac{1}{p^l} ~.
$$
For each $\de > 0$ we define
\beql{Eq4.9}
{\mathcal H}^{s} ( \Om, \de) := \inf \left \{
\sum_{j=1}^{\infty} |I_j|^s ~:~
\Om \sset \bigcup_{j=1}^{\infty} I_j \,,
~ |I_j| \le \de \mbox{~for all~} j \right \} ~,
\eeq
and\footnote{The limit (which may be $\infty$) exists since
$ {\mathcal H}^{s} ( \Om, \de_1) \ge  {\mathcal H}^{s} ( \Om, \de_2) $
if $\de_1 \le \de_2.$}
\beql{Eq4.10}
{\mathcal H}_{p}^{s} ( \Om) := 
\lim_{\de \rightarrow 0} {\mathcal H}^{s} ( \Om, \de) ~.
\eeq
The {\em $p$-adic Hausdorff dimension} of $\Om$ is the unique
value $s_0$ such that ${\mathcal H}_{p}^{s} ( \Om) = \infty$ for $s < s_0$
and ${\mathcal H}_{p}^{s} ( \Om) = 0$ for $s > s_0$.
The value of ${\mathcal H}_{p}^{s_0} ( \Om)$ may be zero, finite or infinite.

A closed set $\Omega \subset \ZZ_p$ is called {\em weakly self-similar}
$(\bmod~p^k)$ 
{\em with branching ratio} $b$, where
$b$ is an integer $\ge 2$,
if the following ``equal branching''  property holds for $l = 1, 2, \ldots$.
Let $W_l( \Omega)$ denote the set of initial sequences
of digits of length $lk$ in $\Omega$, i.e.
\beql{Eq4.10aa}
W_l( \Omega) =       \{
\be = \sum_{j=0}^{lk-1} a_j p^j ~:~
\mbox{~there~ exists~some~} \al \in \Omega 
\mbox{~whose ~``initial part'' }
p^{lk}P_p(p^{-lk}\al) = \beta \} ~.
\eeq
Then each sequence in $W_l(\Omega)$
should extend to exactly $b$ sequences in  $W_{l+1}(\Omega)$.
That is, if the 
 digits of an element $\beta \in \Omega$
are grouped in blocks of size $k$, once the first $l$ blocks of digits are
specified, there are exactly $b$ allowable choices for the next
block of digits. 

\begin{theorem}\label{Th4.1}
For each $k \ge 1$ the set
$p^k \Om_k(p)$   is weakly self-similar (mod $p^k$) 
with branching ratio $b = \phi(p^k) = p^k - p^{k-1}$. 
\end{theorem}

\paragraph{Proof.}
The set $W_l( \Om_k(p))$ specifies the conditions under
which the first $l$ iterates
$f(\al)$, $f^{(2)}(\al)$, $\ldots$, $f^{(l)}(\al)$ $\notin 
\frac{1}{p^{k-1}} \, \ZZ_p$.
Each such condition is a congruence (mod $p^{lk}$),
which has exactly $\phi(p^k) = p^k - p^{k-1}$ solutions for the next digit
(compare \eqn{Eq3.20} in the proof of Theorem \ref{Th1}, taking 
each $d_j = p^k$). Finally the definition of $\Om_k(p)$ implies
it is a closed set. Thus $p^k \Om_k(p)$ is weakly self-similar. 
~~~$\bsq$

We can determine the Hausdorff dimension of weakly self-similar sets 
in $\ZZ_p$, together with upper and lower
bounds for the Hausdorff measure at
this dimension.

\begin{theorem}\label{Th4.2}
Let $\Om \sset   \ZZ_p$ be a weakly
self-similar set (mod $p^k$) with branching ratio $b$
satisfying $2 \le b < p^k$.
Then $\Om$ is a compact set and has Hausdorff dimension
$s(\Omega)$ given by 
$$
\dim_{H} ( \Omega) = \frac{\log b}{\log p^k} ~.
$$
Its $s(\Omega)$-dimensional Hausdorff measure
satisfies
\beql{Eq4.14}
\left ( \frac{b}{p^k} \right) ^ { 1 - \frac{1}{k}}   ~\le~
{\mathcal H}_p^{s(\Omega)} ( \Omega) 
~\le~ 1 ~.
\eeq
\end{theorem}

\paragraph{Proof.}
The compactness of $\Om$ is established similarly
to  Theorem \ref{Th4.1}.
We define $W_l(\Om)$ as in \eqn{Eq4.10aa}. This set has cardinality 
$b^{l}$
by hypothesis, and
$$
\Om = \bigcap_{l=1}^{\infty} \widetilde{W}_l ( \Om) ~,
$$
where  $\widetilde{W}_l ( \Om)$
is the compact set
$$
\widetilde{W}_l ( \Om) = \left \{
\tilde{\be} \in  \ZZ_p ~:~
\tilde{\be} \equiv \be ~(\bmod~ p^{lk}) \mbox{~for~some~} \be \in
W_l(\Om) \right \} ~.
$$

To determine the Hausdorff
dimension it  suffices to establish the inequalities \eqn{Eq4.14},
since the fact that the Hausdorff measure 
is positive and finite determines the Hausdorff dimension.
Set $s = (\log b)/(\log p^k),$ so that
$p^{ks} = b.$ Also $0 < s < 1$.

For the upper bound in \eqn{Eq4.14} we consider the sets
$\widetilde{W}_l ( \Om)$.
Now $b^{l}$ balls of diameter $p^{lk}$
cover $\widetilde{W}_l ( \Om)$, hence cover $\Om$.
Thus for $\de = p^{-lk}$ this covering gives
$$
{\mathcal H}^{s} ( \Om, p^{-lk}) \le  \frac{b^l}{p^{lks}} = \frac{b^l}{b^l} = 1 ~,
$$
which implies 
$$
{\mathcal H}_p^{s} ( \Om))  \le 1 \,.
$$

The lower bound argument is similar in spirit to that 
used for  Cantor sets in
[Falconer 1990, pp. 31--32].
By the compactness of $\Om$ we need only prove that the lower bound
holds for finite coverings.
The non-archimedean property of the valuation $|~|_p$
means that each $I_j$ has diameter $p^m$ for some $m$, and hence we can
enlarge $I_j$ to a ball $B(\al; p^m) \supseteq I_j$
without changing its diameter.  But $B(\al; p^m)$ gives an open cover,
so it has a finite subcover:
$$
\Om \sset \bigcup_{j=1}^{m} B(\al_j; p^{m_j}) ~,
$$
and
$$
\sum_{j=1}^{m} | B(\al_j; p^{m_j}) |^s ~=~ \sum_{j=1}^{m} p^{m_j s} ~.
$$
We want to show
\beql{Eq4.18}
\sum_{j=1}^{m} | B(\al_j; p^{m_j}) |^s ~\ge~ 
\left ( \frac{b}{p^k} \right ) ^{ \frac{k-1}{k}}  ~.
\eeq
We first replace these balls with balls of diameter $p^{-k l_j}$
for integers $l_j$. Write
$$
m_j = -k l_j + k_j \, , ~ 0 \le k_j \le k-1 ~.
$$
Then we can cover $ B(\al_j; p^{m_j})$ with $p^{k_j}$ balls
$B(\al_{j'}; p^{-k l_j})$. We claim that
\beql{Eq4.19}
\sum_{j=1}^{m} | B(\al_j; p^{m_j}) |^s ~\ge~ 
\left ( \frac{b}{p^k} \right ) ^{ \frac{k-1}{k}}
\sum_{j'} | B(\al_{j'}; p^{-k l_j}) |^s ~.
\eeq
This will follow if we show that
$$
\left ( p^{-k l_j + k_j} \right ) ^ s \ge \left (\frac{b}{p^k} \right ) ^ { \frac{k-1}{k}}
p^{k_j} \left ( p ^ {-k l_j } \right ) ^s ~,
$$
holds for each $j$. This in turn is equivalent to showing 
$$
p^{k_j (s-1)} \ge \left ( \frac{b}{p^k} \right ) ^ { \frac{k-1}{k}} ~.
$$
Since $p^{ks}=b$ we have
$$
p^{k_j (s-1)} =
\left ( \frac{b}{p^k} \right ) ^ { \frac{k_j}{k}}
\ge \left ( \frac{b}{p^k} \right ) ^ { \frac{k-1}{k}} ~,
$$
which proves \eqn{Eq4.19}.
Thus \eqn{Eq4.18} will follow from showing
\beql{Eq4.20}
\sum_{j=1}^{m} | B(\al_j; p^{-k l_j}) |^s ~\ge~  1
\eeq
for any set of such balls that covers $\Om$.
We may suppose $l_1 \le l_2 \le \cdots \le l_m$.
By the weak self-similarity of $\Om$
there are $b^{l_m}$ 
principal parts to cover
with balls of diameter $p^{-k l_m}$ in
$\widetilde{W}_m ( \Om)$.
Weak self-similarity also says that each ball of radius
$p^{-k l_j}$ covers either none or else exactly 
$b^{l_m - l_j}$ such
principal parts in $\widetilde{W}_m ( \Om)$.
Since the balls cover $\Om$, we must have
$$
\sum_{j=1}^{m} b^{l_m - l_j} \ge b^{l_m} ~.
$$
Dividing by $b^{l_m}$ yields
\beql{Eq4.21}
\sum_{j=1}^{m} b^{ - l_j} \ge 1 ~.
\eeq
Using this bound we obtain
\begin{eqnarray}
\sum_{j=1}^{m} | B(\al_j; p^{-k l_j}) |^s
& ~=~ & 
\sum_{j=1}^{m} p^{-k l_j s}
\nonumber \\
& ~=~ & 
\sum_{j=1}^{m} b^{ -l_j }
\nonumber \\
& ~\ge~ & 
1  ~,
\nonumber
\end{eqnarray}
which is \eqn{Eq4.20}.
Thus \eqn{Eq4.18} holds.  ~~~$\bsq$

\begin{coro}
The set $\Omega_k(p)$ has Hausdorff dimension
\beql{Eq4.13a}
\dim_{H} ( \Om_k(p)) = s(p^k) ~:=~ 
1 ~-~
\frac { \log \left (  1 + \frac{1}{p-1} \right ) } { k \log p } 
\eeq
for $k \ge 1$.
Furthermore, its $s(p^k)$-dimensional Hausdorff measure
satisfies
\beql{Eq4.14a}
\left ( 1 - \frac{1}{p} \right ) ^ { 1 - \frac{1}{k}}
(p^k - p^{k-1}) ~\le~
{\mathcal H}_p^{s(p^k)} ( \Om_k(p)) 
~\le~ p^k-p^{k-1} ~.
\eeq
\end{coro}

\paragraph{Proof.} This follows by applying Theorems~\ref{Th4.1}
and \ref{Th4.2} to the set $p^k \Omega_k(p)$
and using the fact that
$$ 
{\mathcal H}_p^{s(p^k)} ( \Om_k(p)) = b {\mathcal H}_p^{s(p^k)} (p^k \Om_k(p)) ~,
$$
since $(p^{k})^s =b$.
Using  the branching ratio
$b = \phi(p^k) = p^k - p^{k-1},$
we have
$$
\frac{\log b}{\log p^k} = \frac{\log p^k + \log (1 - 1/p)}{\log p^k}
$$
which gives \eqn{Eq4.13a}. 
~~~$\bsq$

\paragraph{Remark.}
For the branching ratio $b = \phi(p^k) = p^k - p^{k-1}$,
one can show that equality may occur on either side 
of the Hausdorff measure bounds in \eqn{Eq4.14} (or \eqn{Eq4.14a}).
For the lower bound, take the $p^k-p^{k-1}$ allowed
digit sets in each layer to be
$\sum_{j=lk}^{k(l+1)-1} a_j \, p^j \,$
with the restriction that   $a_{lk} \neq 0.$
Then we can cover $\Om$ with
$(p^k-p^{k-1})^{l+1}(p-1)$
balls of radius $p^{-(lk+1)}$, and get
\begin{eqnarray}
{\mathcal H}_p^{s} ( \Om_k(p)) & ~\le~ & (p^k - p^{k-1})(p^k - p^{k-1})^l p^{-lks}
[ (p-1)p^{-s}] \nonumber \\
 & ~\le~ & (p^k - p^{k-1}) \left [ \left ( 1 - \frac{1}{p} \right )
\left ( 1 - \frac{1}{p} \right )^{-\frac{1}{k}} \right ] ~.
\nonumber 
\end{eqnarray}
Since this coincides with the lower bound in \eqn{Eq4.14a},
${\mathcal H}_p^{s} ( \Om_k(p))$ must equal this bound.

For the upper bound, choose the $p^k - p^{k-1}$
allowed digits in each layer to be
$\sum_{j=lk}^{k(l+1)-1} a_j \, p^j$
with the restriction that $a_{k(l+1)-1} \neq 0$.
Then each residue class (mod
$p^{-k(l+1)+r}$) covers exactly $(p-1)p^{r-1}$ classes of $\Om$
(mod $p^{-lk}$), and we can do no better than the upper bound.

It would be interesting to know how many elements
$r = \frac{l}{p^k}$ have no iterate
$f_p^{n}(r) \in \frac{1}{p^{k-1}} \, \ZZ_p$.
We know that this set contains the 
$p^k-p^{k-1}$ elements
$1 \le l \le p^k-1$ with $(l,p) = 1$.
We conjecture that these are the only such elements.

%
%
%
%

\section{Approximate multiplication maps}~\label{Sec5}

We can use similar methods
to study iteration of the approximate multiplication map
$f_r : \QQ \rightarrow \QQ$ given by
\beql{Eq5.1}
f_r(x) = r \lceil x \rceil ~,
\eeq
where $r$ is a fixed rational number, say $r = \frac{l}{d}$
with $\gcd(l, d) = 1$.
In this case, since $x$ enters into the iteration
only as the integer $\lceil x \rceil$, we may restrict
attention to initial values $x \in \ZZ$.
We consider the case
when the denominator $d > 1$,  and study the question
of whether some iterate $f_r^{(j)}(x)$
will be an integer for some $j \ge 1$.
Note that all iterates lie in $\frac{1}{d} \ZZ$.

Unlike the case of approximate squaring, the iterates do
not remain integral once they become integral. However, the
truth of Conjecture~\ref{Con2} would imply that infinitely
many members of a sequence of iterations 
$\{ f_r^{(j)}(n) ~:~ j \ge 1\}$ will be integers, provided $|r| > 1$.

It is convenient to rescale the map to
eliminate the denominators $d$, 
by conjugating $f_r$
by the dilation $\Phi_d(x) = d \, x$. The result is the map
$g_r: \ZZ \to \ZZ$
given by 
$g_r(x) := \Phi_d \circ f_r \circ \Phi_d^{-1}(x)$.
Thus 
\beql{Eq5.1aa}
g_r(x) = d \, f_r(\frac{x}{d}) = l \lceil \frac{x}{d} \rceil ~,
\eeq
with
$$
g_r^{(j)}(x) = d \, f_r^{(j)}(\frac{x}{d})
$$
for $j = 1, 2, \ldots$.
We have 
$$
g_r(n) = \frac{1}{d}( l n +l_b) \mbox{~when~} n \equiv b~ (\bmod~d) ~,
$$ 
where  $l_0 = 0$ and 
$$
l_b = l(d - b)~~ \mbox{for}~~ 1 \le b \le d-1 ~.
$$
For example, when $r= \frac{3}{2}$, we have
\[
g_{3/2}(n) = \left\{
 \begin{array}{ll}
   \frac{3}{2}n &  \mbox{~if~ $n \equiv 0~(\bmod~ 2)$} \, , \\
   \frac{3}{2}n + \frac{3}{2}&  \mbox{~if~ $n \equiv 1~(\bmod~ 2)$}\,.
 \end{array}
\right.
\]
Our question then becomes: when does
the sequence of iterations $\{ g_r^{(j)}(n) ~:~ j \ge 1\}$
contain a term
which is divisible by $d$?

The  map $g_r$ belongs a general class of functions
which we will denote by $\sP_r$, consisting
of those ``periodically linear'' functions $h_r: \ZZ \to \ZZ$
of the form 
\beql{Eq5.1b}
h_r(n) = \frac{1}{d}( l n + l_b) \mbox{~when~} n \equiv b~ (\bmod~d) ~,
\eeq
where $r = l/d$ is rational and the integers 
$\{ l_b:~ 0 \le b \le d-1 \}$ satisfy
the conditions
\beql{Eq5.1c}
l_b \equiv - l b ~(\bmod~d)
\eeq
needed to give an integer-valued map.
Thus, although the notation does not reflect this, $h_r$ is defined
by specifying $r = l/d$ and constants
$l_0, l_1, \ldots, l_{d-1}$ satisfying \eqn{Eq5.1c}. 
We note that any map $h_r$ in $\sP_r$ is ``self-similar'' in the
sense that its linear part $\frac{l}{d} n$ is independent
of the residue class. More general  classes of
periodically linear functions have been
studied in connection with the $3x+1$ problem -- see
\S3.2 of  [Lagarias 1985]; those classes differ from $\sP_r$ in
allowing the linear part of
the map to depend on the residue class $(\bmod~d)$.

Our methods apply generally to the question of whether a particular
function $h_r \in \sP_r$ has an iterate that is divisible
by $d$ (the denominator of $r$).
The behavior of the map $h_r$ depends on
whether $|r| > 1$, the ``expanding map'' case;
$r = \pm 1$, the ``indifferent map'' case;
or $|r| \le 1$, the ``contracting map'' case.
Our motivation comes from the expanding map
case, but the results proved below apply to all cases.

We formulate a general conjecture concerning functions in
this class that are expanding maps, of which 
Conjecture \ref{Con2} is a special case.

\begin{conj}\label{Con4}
Let $r = \frac{l}{d}$ with $\gcd(l, d)=1$ and
$|r| > 1$. Let $h_r: \ZZ \to \ZZ$ be a function in the class
$\sP_r$. Then for each integer
$n$, with at most a finite number of  exceptions,
there is some iterate $j \ge 1$ such that  
$h_r^{(j)}(n) \equiv 0~(\bmod~d)$.
\end{conj}

The ``expanding'' condition on $r$ is necessary, for the
conjecture fails for certain $r$ with $0 < r < 1,$
as shown in Theorem~\ref{th5.3} below.

An interesting map in the class $\sP_r$ is
\beql{Eq5.2}
\tilde{g}_r (x) = \lceil rx \rceil ~.
\eeq
This map has
$l_0 = 0$ and
$$
l_b = d - (lb~\bmod~ d ) \mbox{~for~} 1 \le b \le d-1 ~.
$$
The function $\tilde{g}_{3/2}(x)$ appears in Mahler's study
of $Z$-numbers [Mahler 1968],  as explained below.
We have
\[
\tilde{g}_{3/2}(n) = \left\{
 \begin{array}{ll}
   \frac{3}{2}n &  \mbox{~if~ $n \equiv 0~(\bmod~ 2)$} \, , \\
   \frac{3}{2}n + \frac{1}{2}&  \mbox{~if~ $n \equiv 1~(\bmod~ 2)$}\,.
 \end{array}
\right.
\]
Mahler's study of $Z$-numbers led to
questions similar to Conjecture~\ref{Con4}.
A {\em $Z$-number} is a positive real number
$\xi$ with the property that
$$
0~ \le~ \leftBra \, (\frac{3}{2})^n \xi~ 
\rightBra ~\le~ \frac{1}{2} \mbox{~for~all~}~ n \ge 1 ~,
$$ 
where $\leftBra x  \rightBra = x - \lfloor x \rfloor$ denotes the
fractional part of $x$.
Mahler conjectured that $Z$-numbers do not exist, and showed
that a necessary and sufficient condition for their
non-existence is that for each $n \ge 1$ there exists
some $j \ge 1$ (depending on $n$) such that
\beql{Eq5.2a}
\tilde{g}_{3/2}^{(j)}(n) \equiv 3~(\bmod~4) ~.
\eeq
Mahler's conjecture remains open. Mahler obtained 
a nontrivial  upper bound on the number of $Z$-numbers 
smaller than $x$, and [Flatto 1991] 
improved the upper bound to $O(x^{0.59})$ for $x \to \infty$.

In comparison to \eqn{Eq5.2a},  
Conjecture~\ref{Con4} for 
$ \tilde{g}_{3/2}$ (with $r= \frac{3}{2}$) asserts 
that for each $n \in \ZZ$ there exists some  
$j \ge 1$ (depending on $n$)
with
\beql{Eq5.3}
\tilde{g}_{3/2}^{(j)}(n) \equiv 0~(\bmod~2) ~,
\eeq
with at most a finite number of exceptions.
This special case of Conjecture~\ref{Con4} is true, as a
consequence of
the next theorem. There is exactly one exceptional integer,  $n= -1$, 
whose iterates never satisfy \eqn{Eq5.3}.

More generally, for the  case of rational numbers $r$ 
with denominator $d=2$, Conjecture \ref{Con4}
is provable for all functions in the class $\sP_r$ in a fashion
analogous  to that used for the approximate squaring map
in Section \ref{Sec2}.

\begin{theorem}~\label{th5.1}
Let $h_r$ be a function in the class $\sP_r$,
for fixed  $r = \frac{2 t+ 1}{2}$ where $t$ is an integer. 
Then for
each $n \in \ZZ$, with at most two exceptions, there exists some iterate
$k \ge 1$ with
\beql{Eq5.4}
h_r^{(k)}(n) \equiv 0~(\bmod~ 2) ~.
\eeq
\end{theorem}

\paragraph{Proof.}
We have
$$
h_{r}(n) = \left\{
 \begin{array}{ll}
   r\,n + l_0 &  \mbox{~if~ $n \equiv 0~(\bmod~ 2)$} \,, \\
   r\,n + \frac{1}{2} + l_1 &  \mbox{~if~ $n \equiv 1~(\bmod~ 2)$}\,.
 \end{array}
\right.
$$
for some integers $l_0, l_1$.

We claim that the set of integers $n$ satisfying
\beql{Eq5.5}
h_r^{(j)}(n) \not\equiv 0 ~(\bmod~2) \mbox{~for~} 1 \le j \le k
\eeq
consist of the integers in exactly two arithmetic
progressions $ b~ (\bmod~ 2^{k+1})$, one consisting of
even integers and one of odd integers.
We prove the claim by induction on $k \ge 1$. For the
base case, let $b \equiv a_0~ (\bmod~2)$ with $a_0 = 0$ or $1$ fixed,
and consider the
arithmetic progression $n= b + 2m$, with $m \in \ZZ$.
Then
$$
h_r(n) = h_r(b) + (2t+1)m ~,
$$ 
and the condition
$ 
h_r(n) \equiv 0~ (\bmod~2)
$
then restricts $m$ to lie in a single congruence class 
$m \equiv a_1(b)~ (\bmod~ 2)$. We conclude that exactly
two congruence classes $b \equiv a_0 + 2 a_1~ (\bmod~4)$ satisfy
\eqn{Eq5.5} for $k=1$, one consisting of
even integers and one of odd integers, completing the base case.

For the induction step, supposing \eqn{Eq5.5} true for
$k$, let $b~(\bmod~ 2^{k+1})$
run over the two allowed
congruence classes for the given $k$. Consider the
arithmetic progression $n = b +  2^{k+1} m$,
with $m \in \ZZ$. Then we have
$$
h_r^{(j)} (n) \equiv h_r^{(j)}(b) ~ (\bmod~2)~~\mbox{for}~~ 1 \le j \le k ~,
$$
and
$$
h_r^{(k+1)}(n) = h_r^{k+1}(b) + (2t+1)^{k+1}m ~.
$$
The condition that 
$$
h_r^{(k+1)}(n) \equiv 0 ~(\bmod~2)
$$
is equivalent to  $m \equiv a_{k+1}(b)~ (\bmod~2)$,
which excludes the congruence class
$$b' \equiv b + a_{k+1}(b) 2^{k+1} ~(\bmod~ 2^{k+2}) ~.$$
Thus two congruence classes $(\bmod~ 2^{k+2})$ remain
which satisfy \eqn{Eq5.5} for $1 \le j \le k+1$. 
Since each of the previous  classes $b ~(\bmod~ 2^{k+1})$
contributed one of these classes, one contains even 
integers and the other contains odd integers.
This completes the induction step.

Denote these two classes by $b_0(k+1)~ (\bmod~ 2^{k+2})$ and 
$b_1(k+1)~ (\bmod~ 2^{k+2})$, respectively.
It follows that \eqn{Eq5.4} holds 
except for integers in the sets
$$
\bigcup_{k=1}^{\infty} \{ n \equiv b_0(k)~ (\bmod~ 2^{k+1}) \} ~.
$$
and
$$
\bigcup_{k=1}^{\infty} \{ n \equiv b_1(k)~ (\bmod~ 2^{k+1}) \} ~.
$$
Each of these sets contains at most one element, so there
are at most two exceptional elements.
~~~$\bsq$

Concerning Theorem~\ref{th5.1}, there exist examples of 
functions $h_r \in \sP_r$ with denominator
$d=2$ such that there 
are zero, one or two elements in the exceptional set.
For example, for $r= \frac{3}{2}$
the function $\tilde{g}_{3/2}$ has
the single exceptional element $n = -1$.
The map
\beql{Eq5.4d}
h_{3/2}(n) = \left\{
 \begin{array}{ll}
   \frac{3}{2}\,n + \frac{1}{2} &  \mbox{~if~ $n \equiv 1~(\bmod~ 2)$} \,, \\
   \frac{3}{2}\,n - 1 &  \mbox{~if~ $n \equiv 0~(\bmod~ 2)$} \,
 \end{array}
\right.
\eeq
has two exceptional points, $0$ and $-1$, all the iterates
of which are odd.

We also note that 
Theorem ~\ref{th5.1} applies to the two values  $r = \pm \frac{1}{2}$
where the maps in $\sP_r$ are  contracting.

We next prove a result in the direction of 
Conjecture \ref{Con4} for denominators $d \ge 3$.
We bound the number of exceptions below $x$ for
a general function in the class $\sP_r$, using an argument
similar to that for the approximate
squaring map studied in Section \ref{Sec3}.
Given a rational number $r = \frac{l}{d}$ with $\gcd(l, d) = 1$ 
and a function $h_r$ in the class $\sP_r$, 
we define the exceptional set by
\beql{Eq5.7}
E(h_r; x) := \{ n \in \ZZ:~ |n| \le x,  \,
h_r^{(k)}(n) \not\equiv 0~ (\bmod~d) \mbox{~for~all~} k \ge 1\} ~,
\eeq 
and let
$$
N(h_r; x) := \#( E(h_r; x)) ~.
$$
(The exceptional set in \eqn{Eq1.aa}
is $E(r) = E(g_r; \infty)$.)
The following  result holds for all rational $r$, including
those with $|r| \le 1$.

\begin{theorem}~\label{th5.2}
Let $r= \frac{l}{d}$ be a rational number with $\gcd(l, d) = 1$,
and suppose that $d \ge 2$.

(1) For each function 
$h_r \in \sP_r$, the set of integers $l$ such that 
the function $h_r$ iterated starting from the initial
value $x = l$ reaches an integer for the
first time at iterate $j \ge 1$ is a union of $d(d-1)^{j}$
complete arithmetic progressions $(\bmod~ d^{j+1})$.

(2) There is a constant $0 \le \beta_d < 1$
depending only on $d$ such that  for every function 
$h_r \in \sP_r$ we have
\beql{Eq5.8}
N(h_r; x) \le 4d x^{\beta_d} ~.
\eeq
The precise value of the constant is
\beql{Eq5.9}
\beta_d = \frac{\log (d-1)}{\log d} =
1 - \log_d( 1 + \frac{1}{d-1}) ~.
\eeq
\end{theorem}

\paragraph{Proof.}
The arguments to establish this result 
are simpler than those for
the approximate squaring map because we can use radix
expansions to the fixed base $d$. 
If $n \equiv b~ (\bmod~d)$ then 
\beql{Eq5.12}
h_r(n) = \frac{l}{d} n + \frac{l_b}{d} ~,
\eeq
where  $0 \le l_b \le d-1$ with
$l_b \equiv - l b~ (\bmod~d)$. 

We claim  that the elements $n \in \ZZ$
such that 
\beql{Eq5.13}
h_r^{(j)}(n) \not\equiv 0 ~(\bmod~d) \mbox{~for~} 1 \le j \le k
\eeq
consist of a  certain set of $d(d-1)^{k}$ 
residue classes $(\bmod~d^{k+1})$.
We proceed by induction on $k \ge 1$. For the base case
$k=1$, given $b~ (\bmod~d)$ the elements of the 
arithmetic progression $n = b + dm$ with $m \in \ZZ$ have
$$
h_r(n) = h_r(b) + l m ~.
$$
Since $\gcd(l, m) = 1$,  as $m \in \ZZ$ 
varies these numbers cycle through
every residue class $(\bmod~ d)$. In particular there is
one class $m \equiv a_1(b)~ (\bmod~d)$, say, that gives
$g_r(n) \equiv 0~ (\bmod~d)$. This arithmetic progression
$b + a_1(b) d~ (\bmod~ d^2)$ is ruled out and 
all elements of the remaining
$d-1$ arithmetic progressions $(\bmod~ d^2)$ satisfy \eqn{Eq5.13}
with $k=1$. This completes the base case.

For the induction step, suppose \eqn{Eq5.13} holds for $k$,
and there are $d(d-1)^k$ allowed residue classes $b~(\bmod~d^{k+1})$.
For each of these residue classes consider the
arithmetic progression $n = b + d^{k+1}m$ with $m \in \ZZ$.
Using \eqn{Eq5.12} repeatedly, we have
$$
h_r^{k+1}(n) = h_r^{(k+1)}(b) + l^{k+1}m ~.
$$
Since $\gcd(l, m) =1$ this progression cycles through all residue
classes $(\bmod~d)$, and the
condition $h_r^{k+1}(n) \equiv 0~ (\bmod~d)$ rules out one
residue class $b + a_{k+1}(b)d^{k+1} ~(\bmod~ d^{k+2})$, say. Thus
imposing  the additional condition
$$
h_r^{k+1}(n) \not\equiv 0~ (\bmod~ d)
$$
leaves  $d(d-1)^{k+1}$ allowed residue classes $(\bmod~d^{k+2})$
whose elements satisfy \eqn{Eq5.13} for $k+1$.
This completes the induction step, proving the claim.

The proof of the claim establishes part (1), since the
elements that first have an integer iterate at the $k$-th step 
for $k \ge 1$ were
shown to 
form $d(d-1)^{k-1}$ complete arithmetic progressions $(\bmod~d^{k+1})$. 

To establish part (2),
we will prove a stronger result than stated above. 
Set
\beql{Eq5.10}
E^{\ast}(h_r; x) := \{ n \in \ZZ:~ |n| \le x, \, 
h_r^{(k)}(n) \not\equiv 0~ (\bmod~d) \mbox{~for~} 1 \le k  \le \log_d x\}
\eeq
and let $N^{\ast}(h_r; x) = \#(E^{\ast}(h_r; x)).$
Certainly $N^{\ast}(h_r; x) \ge N(h_r; x)$, so it suffices to establish
\beql{Eq5.11}
N^{\ast}(h_r; x) \le 4d x^{\beta_d} ~,
\eeq
where $\beta_d$ is given in \eqn{Eq5.9}.

Now  suppose $d^k \le x < d^{k+1}$.
We observe that 
applying \eqn{Eq5.13}, with $k$ replaced by $k-1$, shows that
the elements in
$E^{\ast} (h_r; x)$ are necessarily  contained in a set of
$d(d-1)^{k-1}$  residue
classes $(\bmod~ d^k)$, since $\log_d x > k-1$.
It follows that
$$
N^{\ast}(h_r; x) \le d (d-1)^{k-1} \cdot 
2 \bigg\lceil \frac{x}{d^k} \bigg\rceil 
\le 2 d^2(d-1)^{k-1} \le 4d (d-1)^k ~,
$$
since $d \ge 2$. Now $x \ge d^k$, so $k \le \frac{\log x}{\log d}$ and we have
\begin{eqnarray*}
N^{\ast}(h_r; x) &\le &  4d e^{k \log (d-1)}\\
&\le & 4d x^{\frac{\log(d-1)}{\log d}} ~,
\end{eqnarray*}
as asserted. 
$~~~\bsq$

Theorem~\ref{th5.2} has a probabilistic interpretation
analogous to that of Theorem~\ref{Th3.2}. We leave its
formulation to the reader.

We now show that the  upper bound of part (2) of Theorem~\ref{th5.2} is
nearly best possible for certain values of $r$ 
in the interval $0 < r < 1$,
where the map is a contracting map. 
That is, for suitable maps in the class $\sP_r$, 
the upper bound \eqn{Eq5.8}  of Theorem~\ref{th5.2} is 
within a multiplicative constant of the best
possible upper bound. 

\begin{theorem}~\label{th5.3}
Let $r= \frac{1}{d}$ with $d \ge 3$.
Then, for the conjugated approximate multiplication
map $g_r$ of \eqn{Eq5.1aa}, the exceptional set $E(g_r; x)$ has cardinality
\beql{Eq5.14}
N(g_r; x) \ge \frac{1}{d} x^{\beta_d} ~, \mbox{~for~all~} x \ge d ~,
\eeq
with $\beta_d = \frac{\log (d-1)}{\log d}$.
In particular, the full exceptional set $E(g_r; \infty)$
is infinite.
\end{theorem}

\paragraph{Proof.}
Note that for $r = \frac{1}{d}$
the functions $g_r$ and  $\tilde{g}_r$ coincide.
We claim that, for each $k \ge 1$, the subset $\Sigma_k$ of $[1, d^k]$
given by
$$
\{ n: 1 \le n \le d^k, ~n = a_0 + a_1 d +\ldots + a_{k-1} d^{k-1} 
\mbox{~with all~} a_i \not\equiv -1~(\bmod~ d)~ \mbox{~for~} i \ge 1\}~
$$
is contained in $E(g_r; d^k)$.
(The cardinality of  $\Sigma_k$ is
$d(d-1)^{k-1}$.) We prove this by induction on $k \ge 1$.

For the base case $k=1$ we have 
$$\Sigma_1 = \{ 1, 2, \ldots, d-1 \} \subseteq E(g_r; d)
$$ 
because
each $g_r(n) = 1 \not\equiv 0~(\bmod~ d)$ and $1$ is 
a fixed point of $g_r$. Next,
$$
g_r( \Sigma_k) \subset \Sigma_{k-1} ~,
$$
because
$$
g_r(n) = (a_1 + 1) + a_2 d + a_3 d^2 + \cdots + 
a_{k-1} d^{k-2} \in \Sigma_{k-1} ~,
$$
and $1 \le a_1 + 1 \le d-1$ by hypothesis.
Thus $\Sigma_k \subset E(g_r; d^k)$, 
which completes the induction step.

For   $d^k \le x < d^{k+1}$, with $k \ge 1$,  we have 
\beql{Eq5.16}
N(g_r; x) \ge  N(g_r; d^k) \ge \# \Sigma_k = d(d-1)^{k-1} \ge 
\frac{1}{d}x^{\beta_d} ~,
\eeq
as asserted.
$~~~\bsq$

Since the exceptional set $E(r)$
of the approximate multiplication map $f_r$ has
$E(r) = E(g_r; \infty)$ and
$E(g_r, x) \subset E(g_r; \infty)$, 
Theorem~\ref{th5.3}  shows that the conclusion of
Conjecture~\ref{Con2} does not hold for these values of $r$.

It seems plausible that for all values $-1< r <1$
(except $r=0$) there is some  function in the
class  $\sP_r$ for which the conclusion of
Conjecture~\ref{Con4} does not hold. However we do not
attempt to construct such functions here.

We end this section with a result that shows that in
at least one case the approximate multiplication
map based on the floor function,
$$
F_r(x) := r \lfloor x \rfloor ~,
$$
behaves in exactly
the same way as $f_r(x)$.
A version of this result 
was communicated to us by Benoit Cloitre.

\begin{theorem}
Let $r= \frac{d+1}{d}$ with $d \ge 1$.
Then, for any integer $m \ge 1$,
the sequence $F_r^{(j)}(m+d)$, $j \ge 0$,
takes exactly as long to reach an integer as
as the sequence $f_r^{(j)}(m)$, $j \ge 0$.
\end{theorem}

\paragraph{Proof.}
Let $y_j = f_r^{(j)}(m)$ and $Y_j = F_r^{(j)}(m+d)$.
An easy induction argument shows that $Y_j = y_j + d + 1$
for $j \ge 0$.
$~~~\bsq$

%
%
%
%

\section{Numerical results}~\label{Sec6}
The simplest case where we do not know if the approximate
squaring map $f$ of \eqn{Eq1.1} will always reach an integer
is when the starting value $r = l/d$ has denominator
$d=3$.
We wish to determine $\theta(r)$ (say), the smallest
value of $k \ge 0$ for which $f^{(k)}(r)$ is an integer.

Testing any particular value of $r$ is complicated by the fact--already
illustrated in Section \ref{Sec1}--that the iterates grow
so rapidly. This difficulty can be overcome by writing
the $k$-th iterate $l_k/d_k := f^{(k)}(r)$ in ``base d'':
\beql{Eq6.1}
\frac{l_k}{d_k} = 
\sum_{j=-1}^{\infty} a_j(k) d^k ~,
\eeq
where the ``digits'' $a_j(k)$ satisfy $0 \le a_{j}(k) < d$
(compare \eqn{Eq3.13}),
but storing only the terms
in \eqn{Eq6.1} with $j \le M$. That is, we work mod $d^{M+1}$.
As long as $\theta(r) \le M-1$, we get the correct answer
by finding the smallest $k$ for
which $a_{-1}(k) = 0$. If this has not happened by the time
$k$ reaches  $M$ we increase $M$ and repeat.

For denominator $3$ the value $M = 25$ is sufficient to
show that $\theta(l/3)$ is finite
for $3 \le l \le 2000$.
The following table shows what happens for the first few values.
It gives the initial term, the number of steps to reach an integer,
and the integer that is reached.
$$
\begin{array}{rcccccccccc}
\mbox{start}: & \frac{3}{3} & \frac{4}{3} & \frac{5}{3} & \frac{6}{3} & \frac{7}{3} & \frac{8}{3} & \frac{9}{3} & \frac{10}{3} & \frac{11}{3} & \cdots \\
\theta: & 0 & 2 & 6 & 0 & 1 & 1 & 0 & 5 & 2 & \cdots \\
\mbox{reaches}: & 1 & 8 & 1484710602474311520 & 2 & 7 & 8 & 3 & 1484710602474311520 & 220 & \cdots
\end{array}
$$
(These are sequences \htmladdnormallink{A072340}{http://www.research.att.com/cgi-bin/access.cgi/as/njas/sequences/eisA.cgi?Anum=A072340} and \htmladdnormallink{A085276}{http://www.research.att.com/cgi-bin/access.cgi/as/njas/sequences/eisA.cgi?Anum=A085276} in [Sloane 1995--2003].)
In the range $l \le 2000$ large values of $\theta(l/3)$ are scarce.
The first few record values are $\theta(l/3) = 0, 2, 6, 22, 23$,
reached at $l = 3, 4, 5, 28, 1783$ respectively.

Starting values $r = \frac{d+1}{d}$ take longer to
converge--we discussed the cases $d \le 8$ in Section \ref{Sec1}.
The initial values of $\theta((d+1)/d)$ can be found in
sequence \htmladdnormallink{A073524}{http://www.research.att.com/cgi-bin/access.cgi/as/njas/sequences/eisA.cgi?Anum=A073524} in [Sloane 1995--2003].
The first few record values are 
$0, 1, 2, 3, 18, 26, 56, 79, 200, 225, 388, 1444$,
reached at
$d = 1, 2, 3, 4, 5, 11, 19, 31, 37, 67,$ $149,$ $199$ respectively
(sequences \htmladdnormallink{A073529}{http://www.research.att.com/cgi-bin/access.cgi/as/njas/sequences/eisA.cgi?Anum=A073529}, \htmladdnormallink{A073528}{http://www.research.att.com/cgi-bin/access.cgi/as/njas/sequences/eisA.cgi?Anum=A073528}).
R. G. Wilson, v. [Wilson 2002] has checked that $\theta((d+1)/d)$
is finite for $d \le 500$.

It is amusing to note that the record value $1444$ has the following
interpretation: starting with $\frac{200}{199}$ and
repeatedly approximately squaring,
the first integer reached is
$$
200^{2^{1444}} ~,
$$
a number with about $10^{435}$ digits.

The approximate multiplication map $f_r(x)$
of  \eqn{Eq5.1} is easier to compute since 
it grows more slowly.
Let $\theta_r(n)$ denote
the smallest value of $k \ge 1$ for which $f_r^{(k)}(n)$ is an integer.
We give just one example.
This table shows what happens when $f_{4/3}(n)$ is iterated
with starting value $n$:
$$
\begin{array}{rccccccccccccccc}
n: & 0 & 1 & 2 & 3 & 4 & 5 & 6 & 7 & 8 & 9 & 10 & 11 & 12 & \cdots \\
\theta_{4/3}(n): & 1 & 3 & 2 & 1 & 2 & 9 & 1 & 8 & 3 & 1 &  7 &  2 & 1  &  \cdots \\
\mbox{reaches}: & 0 & 4 & 4 & 4 & 8 & 84 & 8 & 84 & 20 & 12 & 84 & 20 & 16 & \cdots
\end{array}
$$
(sequences \htmladdnormallink{A085068}{http://www.research.att.com/cgi-bin/access.cgi/as/njas/sequences/eisA.cgi?Anum=A085068} and \htmladdnormallink{A085071}{http://www.research.att.com/cgi-bin/access.cgi/as/njas/sequences/eisA.cgi?Anum=A085071}).
Large values of $\theta_{4/3}(n)$ are again scarce.
The first few record values are $\theta_{4/3}(n) = 
1, 3, 9, 15, 17, 18, 24, 27, 28, 30, 40$,
reached at
$0, 1, 5, 161, 1772, 3097, 3473,$ $23084,$ $38752,$ $335165, 491729$ respectively
(sequences \htmladdnormallink{A085328}{http://www.research.att.com/cgi-bin/access.cgi/as/njas/sequences/eisA.cgi?Anum=A085328} and \htmladdnormallink{A085330}{http://www.research.att.com/cgi-bin/access.cgi/as/njas/sequences/eisA.cgi?Anum=A085330}).
We thank J. Earls [Earls 2003] for computing the last six
terms in these two
sequences.

All the evidence supports the conjectures made here; it would be nice to know more.

\section*{Acknowledgements}
We thank Roland Bacher, 
Benoit Cloitre, Jason Earls, Tony D. Noe and Robert G. Wilson, v.,
for their help in
contributing and extending several sequences in
[Sloane 1995--2003] during the early stages
of this investigation.
We also thank 
Bela Bajnok,
Daniel Berend
and
Benoit Cloitre
for some useful comments on the manuscript.

\section*{References}
\begin{description}

\item[{[Beardon 1991]}]
A. F. Beardon, {\em Iteration of Rational Functions},
Springer-Verlag, NY, 1991.

\item[{[Choquet 1980]}]
G. Choquet,
Construction effective de suites $(k(3/2)^{n})$.
\'{E}tude des mesures $(3/2)$-stables, 
{\em C. R. Acad. Sci. Paris}, S\'er. A-B {\bf 291}
(1980), A69--A74.

\item[{[Collet and Eckmann 1980]}]
P. Collet and J.-P. Eckmann,
{\em Iterated Maps on the Interval as Dynamical Systems},
Birkh\"{a}user: Boston 1980.

\item[{[Earls 2003]}]
J. Earls,
Comment on sequence \htmladdnormallink{A085328}{http://www.research.att.com/cgi-bin/access.cgi/as/njas/sequences/eisA.cgi?Anum=A085328} 
in [Sloane 1995--2003],
August 14, 2003.

\item[{[Eisele and Hadeler 1990]}]
P. Eisele and K. P. Hadeler,
Game of cards, dynamical systems, and a characterization of the
floor and ceiling functions,
{\em Amer. Math. Monthly}, {\bf 97} (1990), 466--477.

\item[{[Falconer 1990]}]
K. Falconer,
{\em Fractal Geometry:
Mathematical Foundations and Applications}, 
John Wiley and Sons, Chichester, 1990.

\item[{[Federer 1969]}]
H. Federer, {\em Geometric Measure Theory},
Springer-Verlag, NY, 1969.

\item[{[Flatto 1991]}]
L. Flatto, $Z$-numbers and $\beta$-transformations,
in P. Walters, ed., 
{\em Symbolic Dynamics and its Applications (New Haven, CT, 1991)},
Contemp. Math., {\bf 135},
Amer. Math. Soc., Providence, RI, 1992, pp. 181--201.

\item[{[Graham and Yan 1999]}]
R. L. Graham and C. H. Yan,
On the limit of a recurrence relation,
J. Difference Eqn. Appl. {\bf 5} (1999), 71--95.

\item[{[Lagarias 1985]}]
J. C. Lagarias,
\htmladdnormallink{The 3x+1 problem and its generalizations}
{http://www.cecm.sfu.ca/organics/papers/lagarias/paper/html/paper.html},
{\em Amer. Math. Monthly}, {\bf 92} (1985), 3--23.

\item[{[Lagarias 1992]}]
J. C. Lagarias,
Number theory and dynamical systems,
in:
{\em The Unreasonable Effectiveness of Number Theory (S.~A.~Burr, Ed.),}
Proc. Symp. Applied Math.,
No.~46 (1992), 35--72.

\item[{[Mahler 1968]}]
K. Mahler,
An unsolved problem on the powers of $3/2$,
{\em J. Australian Math. Soc.}, {\bf 8} (1968), 313--321.

\item[{[Sloane 1995--2003]}]
N. J. A. Sloane,
{\em \htmladdnormallink{The On-Line Encyclopedia of Integer Sequences}
{http://www.research.att.com/~njas/sequences/}},
published electronically at www.research.att.com/$\sim$njas/sequences/.

\item[{[Tanton 2002]}]
J. S. Tanton,
{\em \htmladdnormallink{A Collection of Research Problems}
{http://www.themathcircle.org/dozen\%20problems.htm}},
published electronically at  \\
www.themathcircle.org/dozen\%20problems.htm.

\item[{[Terras 1976]}]
R. Terras,
A stopping time problem on the positive integers,
{\em Acta Arith.}, {\bf 30} (1976) 241--252.

\item[{[Terras 1979]}]
R. Terras,
On the existence of a density,
{\em Acta Arith.}, {\bf 35} (1979) 101--102.

\item[{[Wilson 2002]}]
R. G. Wilson, v.,
Comment on sequence \htmladdnormallink{A073524}{http://www.research.att.com/cgi-bin/access.cgi/as/njas/sequences/eisA.cgi?Anum=A073524} in [Sloane 1995--2003], September 11, 2002.

\item[{[Wirsching 1998]}]
G. J. Wirsching,
{\em The Dynamical System Generated by the $3n+1$ Function},
Lecture Notes in Math., Vo1. 1681, Springer-Verlag: New York 1998.

\end{description}

\end{document}